\documentclass{amsart}
\usepackage{amssymb,amsmath,amsthm,hyperref}

\usepackage[normalem]{ulem}

\theoremstyle{theorem}
\newtheorem{theorem}{Theorem}
\newtheorem{proposition}[theorem]{Proposition}

\theoremstyle{definition}

\usepackage{cancel,metrix,multirow,bigdelim}
  
  \tikzset{
   rectangle 0/.style = {fill=white},
   rectangle 1/.style = {fill=white},
   rectangle 2/.style = {fill=lightgray},
   rectangle 3/.style = {fill=white},
   rectangle 4/.style = {fill=lightgray},
   rectangle 5/.style = {fill=white},
   rectangle 6/.style = {fill=lightgray},
   rectangle 7/.style = {fill=white},
   rectangle 8/.style = {fill=black!30!brown},
   rectangle 9/.style = {fill=white!30!blue},
   rectangle 10/.style = {fill=orange},
   rectangle 11/.style = {fill=white},
   rectangle 13/.style = {fill=white!80!green},
}

\newcommand\CR[2][]{%
  \begin{tikzpicture}[#1,ultra thick, rounded corners=2pt, scale=0.5]
    \foreach \row [count=\rc] in {#2} {
      \xdef\offset{0} 
      \foreach \col in \row {
         \draw[rectangle \col] (\offset,-\rc) rectangle ++ (\col, -1);
         \xdef\offset{\numexpr\offset+\col\relax}
      }
    }
  \end{tikzpicture}%
}

\newcommand\CRn[2][]{%
  \begin{tikzpicture}[#1,ultra thick, rounded corners=2pt, scale=0.5]
    \foreach \row [count=\rc] in {#2} {
      \xdef\offset{0} 
      \foreach \col in \row {
         \draw[rectangle \col] (\offset,-\rc) rectangle node[inner sep=0pt,font=\footnotesize] {\col} ++ (\col, -1) ;
         \xdef\offset{\numexpr\offset+\col\relax}
      }
    }
  \end{tikzpicture}%
}

\title{Classical Fibonacci Compositions}

\author{Brian Hopkins}
\address{Department of Mathematics and Statistics, 
Saint Peter's University, Jersey City, NJ 07306, USA}
\email{bhopkins@saintpeters.edu}

\begin{document}

\maketitle

\begin{abstract}
There are three long-known types of restricted integer compositions whose counts match the Fibonacci sequence:\ one from ancient India and two from 19th century England.  We give proofs of these enumeration results using tiling arguments and discuss how these can be used in combinatorial proofs of several Fibonacci number identities.  With MacMahon's notion of conjugation, we show that, for every $n \ge 2$, the compositions of $n$ include subsets whose sizes satisfy the Fibonacci recurrence.
\end{abstract}

\section{Introduction}

The Fibonacci sequence is defined by $F_0 = 0$, $F_1 = 1$, and $F_n = F_{n-1}+F_{n-2}$ for $n \ge 2$; the next few terms are $1, 2, 3, 5, 8, 13$.  The numbers come from a problem in the 1202 book \textit{Liber abaci} by Leonardo Pisano, also known as Fibonacci:\ ``How many pairs of rabbits are created by one pair in one year?" \cite{s02}.  Here, we will investigate certain integer compositions whose enumerations match the Fibonacci sequence:\ one from ancient India, which certainly predates Fibonacci, and two from 19th century England (if not earlier).

Given a positive integer $n$, a composition of $n$ is an ordered collection of positive integers, called parts, whose sum is $n$.  For example, the compositions of five are 
\begin{align*}
\{ & (5),  (4,1), (3,2), (3,1,1), (2,3), (2,2,1), (2,1,2), (2,1,1,1), \\ 
& (1,4), (1,3,1), (1,2,2), (1,2,1,1), (1,1,3), (1,1,2,1), (1,1,1,2), (1,1,1,1,1) \}.
\end{align*}
We write $C(5)$ for the set and $c(5) = 16$ for the count $\vert C(5) \vert$.

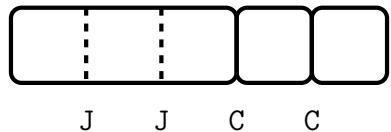
\begin{figure}[h]
\centering
\begin{tikzpicture}[ultra thick, rounded corners=4pt]
\draw[draw=black] (0,0) rectangle ++(3,1);
\draw[draw=black] (3,0) rectangle ++(1,1);
\draw[draw=black] (4,0) rectangle ++(1,1);
\draw [dashed] (1,1) -- (1,0);
\draw [dashed] (2,1) -- (2,0);
\node at (1,-0.5) {\Large \texttt{J}};
\node at (2,-0.5) {\Large \texttt{J}};
\node at (3,-0.5) {\Large \texttt{C}};
\node at (4,-0.5) {\Large \texttt{C}};
\end{tikzpicture}
\caption{The tiling for the composition $(3,1,1)$ and the corresponding cut/join sequence.} \label{311}
\end{figure}

In 1893, Percy MacMahon \cite{m93} used a combinatorial argument to establish a formula for $c(n)$ given $n \ge 1$.  (By convention, $c(0) = 1$.)  Imagine a composition of $n$ as a tiling of a $1 \times n$ board where a part $k$ corresponds to a $1 \times k$ tile.  For instance, Figure \ref{311} shows the tiling for $(3,1,1) \in C(5)$ where the four letters underneath the graphic indicate whether two cells are joined (\texttt{J}) within a single part or cut (\texttt{C}) to separate two parts.  Building a composition of $n$ requires $n-1$ of these binary choices to make a cut/join sequence, and each sequence of $\{\texttt{C}, \texttt{J}\}$ choices leads to a distinct composition, so $c(n) = 2^{n-1}$.  (For more examples, the cut/join sequence for $(2,2,1)$ is \texttt{J C J C} and for $(5)$ is \texttt{J J J J}.  See Figure \ref{conj}, also.) 
 
Rather than allowing all positive integers as parts in compositions, sometimes we restrict which integers can be used as parts.  Such modifications of compositions can be motivated by applications or mathematically interesting results.  The first example of classical Fibonacci compositions comes from an application to literature, when scholars enumerated possible patterns in lines of poetry.  The remaining two classical Fibonacci compositions, due to Augustus De Morgan and Arthur Cayley, have more mathematical inspirations.  After considering the three types of compositions, we discuss their use in combinatorial proofs of identities involving Fibonacci numbers.

\section{Prosody and counting}
Poets work with not just the meanings of words but also structures such as rhyme and rhythm.  Readers of Shakespeare know the rhythmic meter iambic pentameter where each line has ten syllables, specifically, five unstressed--stressed pairs.  For example, Friar Laurence entices Romeo 

\medskip
\begin{center}
\metrics{u p u e p u e p u p u p} {with twen-ty-{-}-hun-dred-{-}-thou-sand times more joy}
\end{center}
(\textit{The Tragedy of Romeo and Juliet}, Act 3, Scene 3), using a baroque expression for two million, where breves indicate unstressed syllables and accents mark stressed syllables.  As the poet Paisley Rekdal explains, ``Prosody is the study of poetic meter, and metrical scansion is what we use to graphically represent changes of sound in a line of metrical verse by marking particular syllables as stressed and others as unstressed'' \cite{r24}.

\begin{figure}[b]
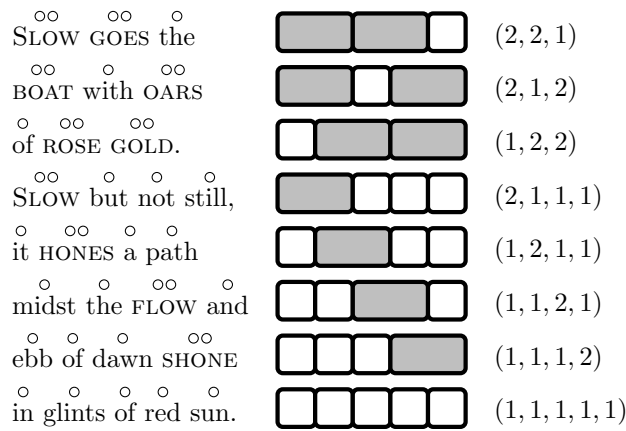

\centering
\renewcommand{\arraystretch}{1.25}
\begin{tabular}{lcl}
\raisebox{3pt}{\metrics{oo oo o}{\textsc{Slow} \textsc{goes} the}} & \CR{{2,2,1}} & \raisebox{4pt}{$(2,2,1)$} \\
\raisebox{3pt}{\metrics{oo o oo}{\textsc{boat} with \textsc{oars}}} & \CR{{2,1,2}} & \raisebox{4pt}{$(2,1,2)$} \\
\raisebox{3pt}{\metrics{o oo oo}{of \textsc{rose} \textsc{gold}}.} & \CR{{1,2,2}} & \raisebox{4pt}{$(1,2,2)$} \\
\raisebox{3pt}{\metrics{oo o o o}{\textsc{Slow} but not still,}} & \CR{{2,1,1,1}} & \raisebox{4pt}{$(2,1,1,1)$} \\
\raisebox{3pt}{\metrics{o oo o o}{it \textsc{hones} a path}} & \CR{{1,2,1,1}} & \raisebox{4pt}{$(1,2,1,1)$} \\
\raisebox{3pt}{\metrics{o o oo o}{midst the \textsc{flow} and}} & \CR{{1,1,2,1}} & \raisebox{4pt}{$(1,1,2,1)$} \\
\raisebox{3pt}{\metrics{o o o oo}{ebb of dawn \textsc{shone}}} & \CR{{1,1,1,2}} & \raisebox{4pt}{$(1,1,1,2)$} \\
\raisebox{3pt}{\metrics{o o o o o}{in glints of red sun.}} & \CR{{1,1,1,1,1}} & \raisebox{4pt}{$(1,1,1,1,1)$}
\end{tabular}
\caption{A poem shown with quantitative meter markings, corresponding tilings of $1 \times 5$ boards, and compositions of five.} \label{qm}
\end{figure}

The poetic analysis relevant to us uses a different metrical system and a nonstandard scansion visualization.  In quantitative meter, ``each syllable is assigned a value that's dependent on its length of sound'' \cite{r24}.  In the poem shown in Figure \ref{qm}, suppose that the long O of the words in small caps lasts twice as long as the other words which all have short vowels.  The syllable lengths are indicated by one or two circles above each word.  Note that each line is five beats long, even though the number of syllables varies from three to five.  Our tiling visualization uses gray $1 \times 2$ dominos for syllables with long vowels and white $1 \times 1$ squares for syllables with short vowels.  The corresponding composition of five is shown on the right.  Notice that the eight lines exhibit every possible tiling of a length five board by dominos and squares, i.e., correspond to the eight compositions of five restricted to parts from $\{1,2\}$.  As we will see, it is not a coincidence that this matches $F_6 = 8$.

Quantitative meter is seldom used in English language poetry, but was applied to verse in Sanskrit and Prakrit, classical languages of South Asia.  Parmanand Singh claims that the writings of Pi\.{n}gala, from perhaps the 5th century BCE, ``indicate a knowledge of the so-called Fibonacci numbers'' \cite{s85}.  Subsequent prosodists Bharata, Virah\={a}\.{n}ka, Gop\={a}la, and Hemacandra displayed increasing knowledge of the recursively defined sequence; all were before Leonardo of Pisa.

Write $C_{12}(n)$ for the compositions of $n$ with parts restricted to $\{1,2\}$, and $c_{12}(n) = \vert C_{12}(n) \vert$.  Figure \ref{qm} includes $C_{12}(5)$ and shows $c_{12}(5) = 8$.

We prove that $c_{12}(n) = F_{n+1}$.  In terms of prosody, this means that there are $F_{n+1}$ different patterns of long and short vowels for a line with $n$ beats.

\begin{proposition} \label{C12}
For $n \ge 1$, the number of compositions of $n$ with parts restricted to $\{1,2\}$ is the $(n+1)$st Fibonacci number, i.e., $c_{12}(n) = F_{n+1}$.
\end{proposition}

\begin{figure}[b]
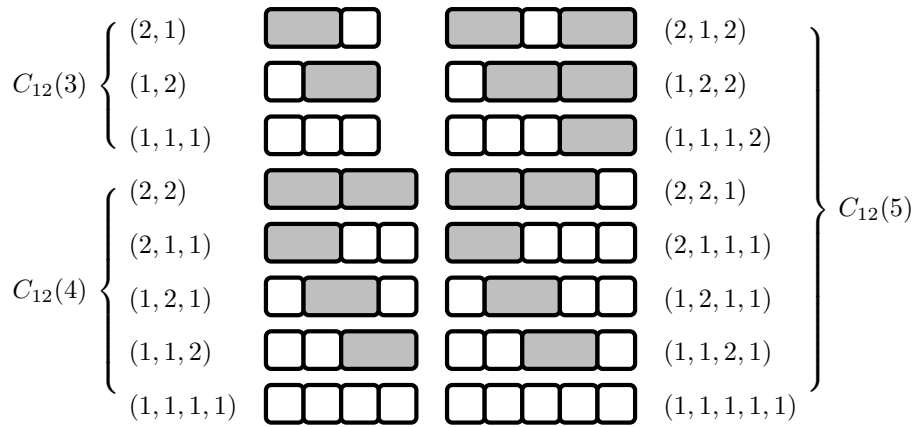

\centering
\renewcommand{\arraystretch}{1.25}
\begin{tabular}{rlllll}
\ldelim\{{3.25}{*}[$C_{12}(3)$ ] \hspace{-1em} & \raisebox{4pt}{$(2,1)$} & \CR{{2,1}} & \CR{{2,1,2}} & \raisebox{4pt}{$(2,1,2)$} & \hspace{-1em} \rdelim\}{9.5}{*}[ $C_{12}(5)$]  \\
& \raisebox{4pt}{$(1,2)$} & \CR{{1,2}} & \CR{{1,2,2}} & \raisebox{4pt}{$(1,2,2)$} \\
& \raisebox{4pt}{$(1,1,1)$} & \CR{{1,1,1}} & \CR{{1,1,1,2}} & \raisebox{4pt}{$(1,1,1,2)$} \\
\ldelim\{{5.5}{*}[$C_{12}(4)$ ] \hspace{-1em} & \raisebox{4pt}{$(2,2)$} & \CR{{2,2}} & \CR{{2,2,1}} & \raisebox{4pt}{$(2,2,1)$} \\
& \raisebox{4pt}{$(2,1,1)$} & \CR{{2,1,1}} & \CR{{2,1,1,1}} & \raisebox{4pt}{$(2,1,1,1)$} \\
& \raisebox{4pt}{$(1,2,1)$} & \CR{{1,2,1}} & \CR{{1,2,1,1}} & \raisebox{4pt}{$(1,2,1,1)$} \\
& \raisebox{4pt}{$(1,1,2)$} & \CR{{1,1,2}} & \CR{{1,1,2,1}} & \raisebox{4pt}{$(1,1,2,1)$} \\
& \raisebox{4pt}{$(1,1,1,1)$} & \CR{{1,1,1,1}} & \CR{{1,1,1,1,1}} & \raisebox{4pt}{$(1,1,1,1,1)$} 
\end{tabular}
\caption{The bijection of Proposition \ref{C12} for $n = 5$.} \label{C12ex}
\end{figure}

\begin{proof}
We establish a bijection $C_{12}(n) \cong C_{12}(n-1) \cup C_{12}(n-2)$ for $n \ge 3$.  That implies $c_{12}(n) = c_{12}(n-1) + c_{12}(n-2)$ and the result will follow from the initial values $c_{12}(1) = 1 = F_2$ from the composition $(1)$ and $c_{12}(2) = 2 = F_3$ from the compositions $(2)$ and $(1,1)$.

To transform a composition in $C_{12}(n-2)$ into a composition of $n$, add an additional part 2 at the end.  In terms of the tiling, add a domino at the end.  

Similarly, for a composition in $C_{12}(n-1)$, add an additional part 1 at the end; equivalently, add a square at the end of the tiling. 

The resulting compositions are all in $C_{12}(n)$ and are distinct:\ Compositions built from $C_{12}(n-2)$ all have final part 2 while compositions from $C_{12}(n-1)$ all have final part 1.  Therefore, \[C_{12}(n-1) \cup C_{12}(n-2) \subseteq C_{12}(n).\]

For the reverse map, for each composition in $C_{12}(n)$, remove its final part/last tile.  For compositions ending in 1 (last tile a square), the images are in $C_{12}(n-1)$, while for compositions ending in 2 (last tile a domino), the images are in $C_{12}(n-2)$.  Therefore, \[C_{12}(n) \subseteq C_{12}(n-1) \cup C_{12}(n-2).\]

The maps are clearly inverses and establish $C_{12}(n) \cong C_{12}(n-1) \cup C_{12}(n-2)$.
\end{proof}

See Figure \ref{C12ex} for an example of the bijection.

The popular book of Benjamin and Quinn \cite{bq03} makes thorough use of square and domino tilings for combinatorial proofs of an impressive array of Fibonacci number identities.

\section{Victorian compositions}

Augustus De Morgan's textbook \textit{Elements of Arithmetic} went through many editions and had a significant impact on British mathematics education.  The fifth edition of 1846 \cite{d46} increased the page count by almost forty percent with appendices on topics ranging from book-keeping and ``casting out nines'' to Horner's method.  The section ``On combinations'' includes the problem
\begin{quote}
Required the number of ways in which a number can be compounded of odd numbers, different orders counting as different ways.
\end{quote}

Write $C_o(n)$ for the compositions of $n$ with all parts odd, and $c_o(n) = \vert C_o(n) \vert$.  Examining $C(5)$ we find
$ C_o(5) = \{(5), (3,1,1), (1,3,1), (1,1,3), (1,1,1,1,1) \}$
so that $c_o(5) = 5$.  Recall that $F_5 = 5$\dots

\begin{proposition} \label{Codd}
For $n \ge 1$, the number of compositions of $n$ with only odd parts is the $n$th Fibonacci number, i.e., $c_o(n) = F_n$.
\end{proposition}

Note that, unlike the situation for $C_{12}(n)$, we cannot add a part 2 to a composition in $C_o(n-2)$ to produce a composition in $C_o(n)$.  

After the thorough proof of Proposition \ref{C12}, the following proofs omit some details.

\begin{figure}[b]
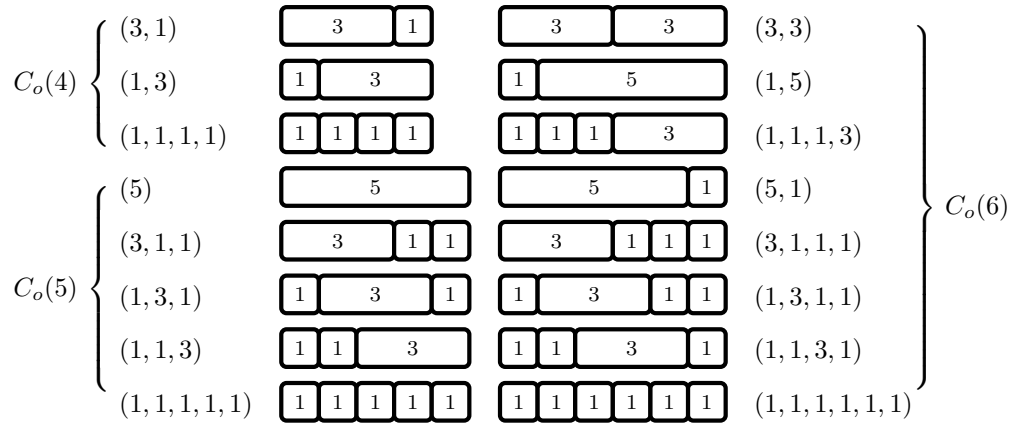

\centering
\renewcommand{\arraystretch}{1.25}
\begin{tabular}{rlllll}
\ldelim\{{3.25}{*}[$C_o(4)$ ] \hspace{-1em} & \raisebox{4pt}{$(3,1)$} & \CRn{{3,1}} & \CRn{{3,3}} & \raisebox{4pt}{$(3,3)$} & \hspace{-1em}\rdelim\}{9.5}{*}[ $C_o(6)$] \\
& \raisebox{4pt}{$(1,3)$} & \CRn{{1,3}} & \CRn{{1,5}} & \raisebox{4pt}{$(1,5)$} \\
& \raisebox{4pt}{$(1,1,1,1)$} & \CRn{{1,1,1,1}} & \CRn{{1,1,1,3}} & \raisebox{4pt}{$(1,1,1,3)$} \\
\ldelim\{{5.5}{*}[$C_o(5)$ ] \hspace{-1em} & \raisebox{4pt}{$(5)$} & \CRn{{5}} & \CRn{{5,1}} & \raisebox{4pt}{$(5,1)$} \\
& \raisebox{4pt}{$(3,1,1)$} & \CRn{{3,1,1}} & \CRn{{3,1,1,1}} & \raisebox{4pt}{$(3,1,1,1)$} \\
& \raisebox{4pt}{$(1,3,1)$} & \CRn{{1,3,1}} & \CRn{{1,3,1,1}} & \raisebox{4pt}{$(1,3,1,1)$} \\
& \raisebox{4pt}{$(1,1,3)$} & \CRn{{1,1,3}} & \CRn{{1,1,3,1}} & \raisebox{4pt}{$(1,1,3,1)$} \\
& \raisebox{4pt}{$(1,1,1,1,1)$} & \CRn{{1,1,1,1,1}} & \CRn{{1,1,1,1,1,1}} & \raisebox{4pt}{$(1,1,1,1,1,1)$} 
\end{tabular}
\caption{The bijection of Proposition \ref{Codd} for $n = 6$.} \label{Coddex}
\end{figure}

\begin{figure}[b]
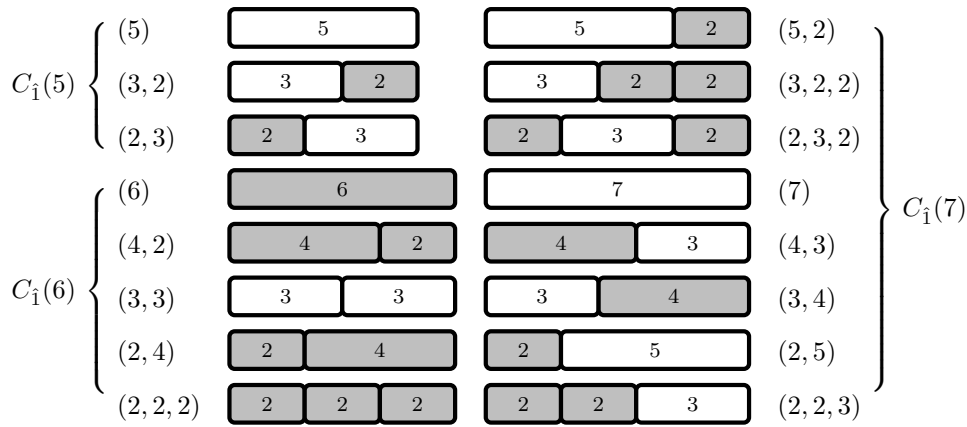

\centering
\renewcommand{\arraystretch}{1.25}
\begin{tabular}{rlllll}
\ldelim\{{3.25}{*}[$C_{\hat{1}}(5)$ ] \hspace{-1em} & \raisebox{4pt}{$(5)$} & \CRn{{5}} & \CRn{{5,2}} & \raisebox{4pt}{$(5,2)$} & \hspace{-1em} \rdelim\}{9.5}{*}[ $C_{\hat{1}}(7)$]\\
& \raisebox{4pt}{$(3,2)$} & \CRn{{3,2}} & \CRn{{3,2,2}} & \raisebox{4pt}{$(3,2,2)$} \\
& \raisebox{4pt}{$(2,3)$} & \CRn{{2,3}} & \CRn{{2,3,2}} & \raisebox{4pt}{$(2,3,2)$} \\
\ldelim\{{5.5}{*}[$C_{\hat{1}}(6)$ ] \hspace{-1em} & \raisebox{4pt}{$(6)$} & \CRn{{6}} & \CRn{{7}} & \raisebox{4pt}{$(7)$} \\
& \raisebox{4pt}{$(4,2)$} & \CRn{{4,2}} & \CRn{{4,3}} & \raisebox{4pt}{$(4,3)$} \\
& \raisebox{4pt}{$(3,3)$} & \CRn{{3,3}} & \CRn{{3,4}} & \raisebox{4pt}{$(3,4)$} \\
& \raisebox{4pt}{$(2,4)$} & \CRn{{2,4}} & \CRn{{2,5}} & \raisebox{4pt}{$(2,5)$} \\
& \raisebox{4pt}{$(2,2,2)$} & \CRn{{2,2,2}} & \CRn{{2,2,3}} & \raisebox{4pt}{$(2,2,3)$} 
\end{tabular}
\caption{The bijection of Proposition \ref{Chat1} for $n = 7$.} \label{Chat1ex}
\end{figure}

\begin{proof}
We establish a bijection $C_o(n) \cong C_o(n-1) \cup C_o(n-2)$ for $n \ge 3$.  The result will follow from the initial values $c_o(1) = 1 = F_1$ from the composition $(1)$ and $c_o(2) = 1 = F_2$ from the composition $(1,1)$.

To transform a composition in $C_o(n-2)$ into a composition in $C_o(n)$, change the last part $k$ to $k+2$; since $k$ is odd, so is $k+2$.  In terms of the tiling, this extends the last tile to be two units longer.  

For a composition in $C_o(n-1)$, add an additional part 1 at the end; equivalently, add a square at the end of the tiling.  

The resulting compositions are all in $C_o(n)$ and are distinct:\ Compositions built from $C_o(n-2)$ all have final part at least three while compositions from $C_o(n-1)$ all have final part 1.

For the reverse map, given a composition in $C_o(n)$ with last part 1 (last tile a square), remove it to produce a composition in $C_o(n-1)$.  Given a composition in $C_o(n)$ with last part $k \ge 3$, replace it with $k-2$ (reduce the last tile's length, which is at least three, by two) to produce a composition in $C_o(n-2)$.
\end{proof}

See Figure \ref{Coddex} for an example of the bijection.  We have included a numeral for each tile's length.

In a short 1876 note, Arthur Cayley considered ``the number of partitions of $n$, no part less than 2, order attended to'' \cite{c76}.  (Now, partitions refers to a related concept where the order of parts is not ``attended to.''  The language distinguishing compositions and partitions was set by MacMahon.)

Write $C_{\hat{1}}(n)$ for the compositions of $n$ with all parts at least two, and $c_{\hat{1}}(n) = \vert C_{\hat{1}}(n) \vert$.  Examining $C(5)$ we find
$ C_{\hat{1}}(5) = \{(5), (3,2), (2,3) \}$
so that $c_{\hat{1}}(5) = 3$.

\begin{proposition} \label{Chat1}
For $n \ge 2$, the number of compositions of $n$ with parts at least two is the $(n-1)$st Fibonacci number, i.e., $c_{\hat{1}}(n) = F_{n-1}$.
\end{proposition}

Note that, unlike the situation for $C_{12}(n)$ and $C_o(n)$, we cannot add a part 1 to a composition in $C_{\hat{1}}(n-1)$ to produce a composition in $C_{\hat{1}}(n)$.

\begin{proof}
We establish a bijection $C_{\hat{1}}(n) \cong C_{\hat{1}}(n-1) \cup C_{\hat{1}}(n-2)$ for $n \ge 4$.  The result will follow from the initial values $c_{\hat{1}}(2) = 1 = F_1$ from the composition $(2)$ and $c_{\hat{1}}(3) = 1 = F_2$ from the composition $(3)$.

To transform a composition in $C_{\hat{1}}(n-2)$ into a composition in $C_{\hat{1}}(n)$, add an additional part 2 at the end.  In terms of the tiling, add a domino at the end.

For a composition in $C_{\hat{1}}(n-1)$, change the last part $k$ to $k+1$; equivalently, extend the last tile to be one unit longer.  

The resulting compositions are all in $C_{\hat{1}}(n)$ and are distinct:\ Compositions built from $C_{\hat{1}}(n-2)$ all have final part 2 while compositions from $C_{\hat{1}}(n-1)$ all have final part at least three.

For the reverse map, given a composition in $C_{\hat{1}}(n)$ with last part 2 (last tile a domino), remove it to produce a composition in $C_{\hat{1}}(n-2)$.  Given a composition in $C_{\hat{1}}(n)$ with last part $k \ge 3$, replace it with $k-1$ (reduce the length of the final tile, which is longer than a domino, by one) to produce a composition in $C_{\hat{1}}(n-1)$.
\end{proof}

See Figure \ref{Chat1ex} for an example of the bijection.  In addition to the numerals indicating tile length, even length tiles are shaded gray.

One can imagine that the results of Propositions \ref{Codd} and \ref{Chat1} were known before De Morgan and Cayley, respectively.  Searching for earlier results is complicated by not only the various terminology for the objects being counted, but also the fact that the name ``Fibonacci numbers'' seems to have first been used in the late 1800s.

\section{Combinatorial proofs}

One motivation for finding different compositions counted by the Fibonacci numbers is providing multiple tools for combinatorial proofs.  For example,
\begin{equation}
F_1 + F_3 + \cdots + F_{2n-1} = F_{2n} \label{fodd}
\end{equation}
has straightforward combinatorial proofs for all three classical Fibonacci compositions as detailed in the author's book \cite[pp.\ 172--173]{h25}.

Some identities seem to have more direct proofs for a particular choice of Fibonacci composition.  We believe that the related identity
\begin{equation}
F_2 + F_4 + \cdots + F_{2n} = F_{2n+1} - 1 \label{feven}
\end{equation}
is most easily shown with $C_o(n)$, De Morgan's compositions \cite[p.\ 168--169]{h25} while
\begin{equation}
F_1 + F_4 + \cdots + F_{3n-2} = \tfrac{1}{2}F_{3n}, \label{fthrover}
\end{equation}
with all terms doubled, can be shown directly with $C_{\hat{1}}(n)$, Cayley's compositions  \cite[p.\ 169]{h25}.
These are all matters of taste, of course, but we believe that the combinatorial proof for \textit{The Book} of Paul Erd\H{o}s's dreams \cite{az18} for the remarkable identity
\begin{equation}
\binom{n-1}{0} + \binom{n-2}{1} + \binom{n-3}{2} + \cdots = F_n, \label{PasFib}
\end{equation}
with a suitable interpretation of binomial coefficients from MacMahon's combinatorial method discussed in the first section, also comes from compositions with parts all at least two \cite[pp. 162--164, 174--175]{h25}. 

Benjamin and Quinn \cite{bq03} prove all of these with their system equivalent to $C_{12}(n)$ (although they express results in terms of $f_n = F_{n+1}$):\ \eqref{fodd} is their Identity 2, \eqref{feven} their Identity 12, \eqref{fthrover} their Identity 23, and \eqref{PasFib} their Identity 4.  It is interesting to work out combinatorial proofs using De Morgan's or Cayley's compositions for their fifty or so identities incorporating Fibonacci numbers.

\begin{figure}
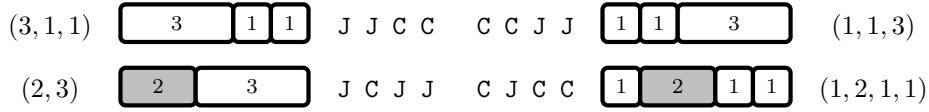

\centering
\renewcommand{\arraystretch}{2}
\begin{tabular}{cccccc}
\raisebox{4pt}{$(3,1,1)$} & \CRn{{3,1,1}} & \raisebox{4pt}{\texttt{J J C C}} & \; \raisebox{4pt}{\texttt{C C J J}} & \CRn{{1,1,3}} & \raisebox{4pt}{$(1,1,3)$} \\
\raisebox{4pt}{$(2,3)$} & \CRn{{2,3}} & \raisebox{4pt}{\texttt{J C J J}} & \; \raisebox{4pt}{\texttt{C J C C}} & \CRn{{1,2,1,1}} & \raisebox{4pt}{$(1,2,1,1)$}
\end{tabular}
\caption{The compositions $(3,1,1)$ and $(2,3)$, their conjugates, and the cut/join sequences.} \label{conj}
\end{figure}

Our final result uses an additional concept from MacMahon \cite{m93}.  Given the cut/join sequence of a composition, a natural operation is to swap each binary choice.  This produces what he calls the conjugate composition; see Figure \ref{conj} for examples.  (Note that the conjugate of $(3,1,1)$ matches reversing the order of $(3,1,1)$.  MacMahon uses the term inverse conjugates to describe pairs such as $(3,1,1)$ and $(1,1,3)$.)

We conclude with another combinatorial result:\ For each integer $n \ge 2$, the compositions of $n$ include an instance of the Fibonacci recurrence $F_{n+1} = F_n + F_{n-1}$.  

\begin{theorem} \label{newthm}
For each $n \ge 2$, there is a bijection $C_{12}(n) \cong C_{\hat{1}}(n) \cup C_o(n)$.
\end{theorem}

Of course, the related numerical result
\[ c_{12}(n) = F_{n+1} = F_n + F_{n-1} = c_o(n) + c_{\hat{1}}(n) \]
follows from Propositions \ref{C12}, \ref{Codd}, and \ref{Chat1}.  What may be new here is the bijection.

Note that $C_{\hat{1}}(n)$ and $C_o(n)$ are not necessarily disjoint sets of compositions.

\begin{figure}[b]
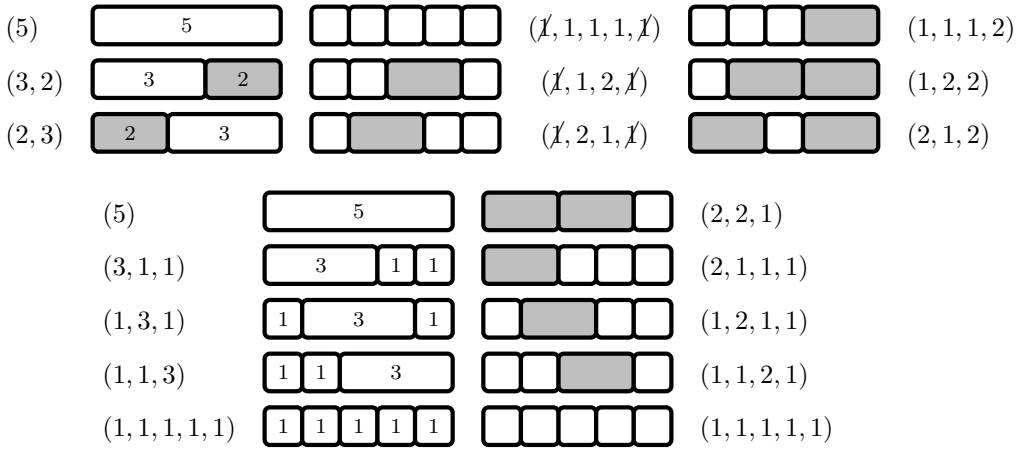

\centering
\renewcommand{\arraystretch}{1.25}
\begin{tabular}{lccccl}
\raisebox{4pt}{$(5)$} & \CRn{{5}} & \CR{{1,1,1,1,1}} &  \raisebox{4pt}{$(\cancel{1},1,1,1,\cancel{1})$}  & \CR{{1,1,1,2}} & \raisebox{4pt}{$(1,1,1,2)$} \\
\raisebox{4pt}{$(3,2)$} & \CRn{{3,2}} & \CR{{1,1,2,1}} & \raisebox{4pt}{$(\cancel{1},1,2,\cancel{1})$} & \CR{{1,2,2}} & \raisebox{4pt}{$(1,2,2)$} \\
\raisebox{4pt}{$(2,3)$} & \CRn{{2,3}} & \CR{{1,2,1,1}} & \raisebox{4pt}{$(\cancel{1},2,1,\cancel{1})$} & \CR{{2,1,2}} & \raisebox{4pt}{$(2,1,2)$}
\end{tabular}

\bigskip
\renewcommand{\arraystretch}{1.25}
\begin{tabular}{lccl}
\raisebox{4pt}{$(5)$} & \CRn{{5}} & \CR{{2,2,1}} & \raisebox{4pt}{$(2,2,1)$} \\
\raisebox{4pt}{$(3,1,1)$} & \CRn{{3,1,1}} & \CR{{2,1,1,1}} & \raisebox{4pt}{$(2,1,1,1)$} \\
\raisebox{4pt}{$(1,3,1)$} & \CRn{{1,3,1}} & \CR{{1,2,1,1}} & \raisebox{4pt}{$(1,2,1,1)$} \\
\raisebox{4pt}{$(1,1,3)$} & \CRn{{1,1,3}} & \CR{{1,1,2,1}} & \raisebox{4pt}{$(1,1,2,1)$} \\
\raisebox{4pt}{$(1,1,1,1,1)$} & \CRn{{1,1,1,1,1}} & \CR{{1,1,1,1,1}} & \raisebox{4pt}{$(1,1,1,1,1)$} 
\end{tabular}
\caption{The bijection of Theorem \ref{newthm} for $n = 5$ with $C_{\hat{1}}(5)$ above and $C_o(5)$ below.} \label{newthmex}
\end{figure}

\begin{proof}
Given $c \in C_{\hat{1}}(n)$, take its conjugate $c'$.  We claim that $c^\prime \in C_{12}(n)$:\ Since each part of $c$ is at least two, there are no adjacent \texttt{C} terms in its cut/join sequence, so the cut/join sequence of $c^\prime$ has no adjacent \texttt{J} terms, i.e., no part is greater than two.  Further, the first and last parts of $c^\prime$ are 1, since $c$ begins and ends with parts at least two.  We send $c$ to the following composition in $C_{12}(n)$:\ Delete the first and last parts 1 of $c^\prime$ and add a part 2 at the end.

Given $c \in C_o(n)$, replace each part $2k+1$ with a sequence of $k$ parts 2 followed by a part 1.  This gives a composition in $C_{12}(n)$ with final part 1.

Considering the last part shows that the images of $C_{\hat{1}}(n)$ and $C_o(n)$ in $C_{12}(n)$ are distinct (even though $C_{\hat{1}}(n) \cap C_o(n)$ may be nonempty).

The reverse map undoes the operations described above.  

For a composition in $C_{12}(n)$ with last part 2, remove that 2, add parts 1 at the beginning and end, then take the conjugate to produce a composition in $C_{\hat{1}}(n)$.  

For a composition in $C_{12}(n)$ with last part 1, convert each run of $k$ parts 2 followed by a part 1 (for $k \ge 0$) to $2k+1$ (since the composition has last part 1, every run of parts 2 is followed by a part 1); this produces a composition in $C_o(n)$.
\end{proof}

See Figure \ref{newthmex} for an example of the bijection.

\end{document}